\documentclass[12pt]{article}

\usepackage{algorithm}
\usepackage{algpseudocode}
\usepackage[hidelinks]{hyperref}
\usepackage{graphicx}

\newtheorem{theorem}{Theorem}

\newtheorem{definition}[theorem]{Definition}

\newtheorem{conjecture}[theorem]{Conjecture}

\begin{document}

	\title{On The Asymptotic Density Of Prime $k$-tuples and a Conjecture of Hardy and Littlewood}
	\author{L\'aszl\'o T\'oth \\ Rue des Tanneurs 7, \\ L-6790 Grevenmacher, G. D. Luxembourg \\ uk.laszlo.toth@gmail.com}
	\date{\textit{Keywords}: Prime $k$-tuple, asymptotic density, conjecture, Skewes number}

	\maketitle
	
	\begin{abstract}
		In 1922 Hardy and Littlewood proposed a conjecture on the asymptotic density of admissible prime $k$-tuples. In 2011 Wolf computed the ''Skewes number'' for twin primes, i.e., the first prime at which a reversal of the Hardy-Littlewood inequality occurs. In this paper, we find ''Skewes numbers'' for $8$ more prime $k$-tuples and provide numerical data in support of the Hardy-Littlewood conjecture. Moreover, we present several algorithms to compute such numbers.
	\end{abstract}
	
	\section{Introduction} \label{intro}
	
	{\let\thefootnote\relax\footnote{{\textit{2010 Mathematics Subject Classification}: Primary 11N05; Secondary 11-04.}}}
	
	Let $(a_1,a_2,\ldots,a_k)$ denote a monotonically increasing sequence of positive even integers and let $p$ be a prime number. Then if the numbers $p+a_i$ for all $1\leq i\leq k$ are prime, the sequence $P=(p,p+a_1,p+a_2,\ldots,p+a_k)$ is a \textit{prime $k$-tuple}. Moreover, if these numbers do not form a complete residue class with respect to any prime, $P$ is \textit{admissible}. Hardy and Littlewood \cite{HardyLittlewood22} made several conjectures concerning the infinitude of admissible prime $k$-tuples. In particular, they conjectured that their asymptotic density can be calculated in terms of the $(a_1,a_2,\ldots,a_k)$ as follows.
	\begin{conjecture} [Asymptotic density of prime $k$-tuples] \label{conj-asym-dens}
		Let $P=(p,p+a_1,p+a_2,\ldots,p+a_k)$ denote an admissible prime $k$-tuple and let $\pi_P(n)$ denote the number of primes $p$ less then a positive integer $n$ such that for all $1\leq i\leq k$, $p+a_i$ is prime. Then
		$$
		\pi_P(n) \sim C_{a_1,a_2,...,a_k}\int_2^n \frac{dt}{\log^{k+1}t},
		$$
		where $C_{a_1,a_2,...,a_k}$ is a constant obtained through a product over all primes $q$ greater than $2$ and the amount of distinct residues of $a_1,a_2,...,a_k$ modulo $q$, denoted by $w(q;a_1,a_2,...,a_k)$, as follows:
		$$
		C_{a_1,a_2,...,a_k}=2^k\prod_q \frac{1- \frac{w(q;a_1,a_2,...,a_k)}{q}}{(1-\frac 1 q)^{k+1}}.
		$$
	\end{conjecture}
	For example, when $k=1$ and $a_1=2$ (i.e., the twin primes), we have $C_2=2\prod_{q}\left(1-\frac{1}{(q-1)^2}\right)= 1.320323632\ldots$, called the twin prime constant (although some authors prefer defining $\frac{C_2}{2}$ as the twin prime constant).
	
	If the density described in Conjecture \ref{conj-asym-dens} is true, it immediately implies the infinitude of the corresponding prime tuple, so any numerical data in its favor is of value. There are several ways to study this density, for instance by analyzing the amount of sign changes in the difference 
	\begin{equation}
	\delta_P(n)=\pi_P(n)-C_{a_1,a_2,...,a_k}\int_2^n \frac{dt}{\log^{k+1}t}.
	\end{equation}
	This is comparable to the study of the so-called Skewes number, the first positive integer $n$ such that $\pi(n)>$Li$(n)$, where Li$(n)$ denotes the logarithmic integral $\displaystyle \int_{2}^{n}\frac{dt}{\log t}$ and $\pi(n)$ is the usual prime-counting function. The existence of such a number was shown first by Littlewood \cite{Littlewood14} and an upper bound was given by Skewes (\cite{Skewes33} assuming the Riemann Hypothesis, \cite{Skewes55} without assuming it). In this paper, we extend this definition by formally defining the ''Skewes number'' for prime $k$-tuples as follows.
	
	\begin{definition} [Skewes number for prime $k$-tuples]
		Let $p$ denote a prime number, $k$ a positive nonzero integer and let $a_i$ for all $1\leq i\leq k$ be positive integers. Furthermore, let $P=(p,p+a_1,p+a_2,\ldots,p+a_k)$ denote an admissible prime $k$-tuple, $\pi_P(n)$ the counting function for $P$ and $C_{a_1,a_2,...,a_k}$ the associated prime tuple constant. Then the first prime $p\in P$ that violates
		$$
		\displaystyle \pi_P(n)<C_{a_1,a_2,...,a_k}\int_2^n \frac{dt}{\log^{k+1}t},
		$$
		if such a prime exists, is the \emph{Skewes number} for $P$.
	\end{definition}
	
	The Skewes number for twin primes was already studied by several authors, among whom we cite Brent \cite{Brent75} and Wolf \cite{Wolf11} as particularly interesting. First, define $P_2=(p,p+2)$ and let $\pi_2(n)$ denote the amount of primes not exceeding $n$ such that $n+2$ is also prime. Furthermore, define Li$\displaystyle_2(n)=\int_{2}^{n}\frac{dt}{\log^2 t}$. Wolf \cite{Wolf11} computed the amount of sign changes in $\delta_{P_2}=\pi_2(n)-C_2$Li$_2(n)$ and found that there are exactly $477118$ up to $2^{48}\approx2.81\times10^{14}$. Moreover, the author found the Skewes number for twin primes, $1369391$. Based on numerical data, Wolf then proposed the following conjecture.
	\begin{conjecture} [Wolf] \label{conj-wolf}
		The amount of sign changes in $d_2(n)=\pi_2(n) - C_2 \rm Li$$_2(n)$ within the interval $n\in(1,T)$ is given by $\frac{\sqrt T}{\log T}$.
	\end{conjecture}
	Note that this result was obtained after a number of ''trials'' and no heuristic argument is given in support of the conjecture. The author's results continue the work of Brent \cite{Brent75}, who also studied the difference $\delta_{P_2}(n)$ and examined its maxima and minima in various intervals $[a,b]$ up to $8\times10^{10}$ by defining the functions:
	$$
	R_3(a,b) = \max\limits_{p \in P_2 \cap [a,b]} \delta_{P_2}(p)
	$$
	and
	$$
	\rho_3(a,b) = \min\limits_{p \in P_2 \cap [a,b]} \delta_{P_2}(p)
	$$
	(note that we adapted Brent's notation to ours in the above definitions). The author also used the numerical data he obtained in order to give an approximation to Brun's constant. Nicely \cite{Nicely04} considered three additional $k$-tuples: $(p,p+2,p+6)$, $(p,p+4,p+6)$ and $(p,p+2,p+6,p+8)$, and computed the corresponding coefficients in their conjectured asymptotic density (denoted by $C$ in Conjecture \ref{conj-asym-dens}) to a high precision. He did not compute the Skewes number for these tuples. Moreover, none of the authors mentioned above provided computer code.
	
	In another paper, motivated by the fact Conjecture \ref{conj-asym-dens} implies that some prime tuples have equal asymptotic density, Wolf \cite{Wolf98} investigated the asymptotic densities of the twin and cousin primes $(p,p+4)$. Indeed, we expect that 
	\begin{equation}
	\pi_{twin}(n),\pi_{cousin}(n) \sim \displaystyle 2 \prod_{p\geq3} \frac{p(p-2)}{(p-1)^2} \int_{2}^{n}\frac{dt}{\log^2 t}.
	\end{equation}
	Wolf looked at the relationship between these two densities by defining the function $y(x)$ as the difference between the number of twin and cousin primes up to a positive integer $x$. Along with a number of numerical computations, Wolf studied the fractal properties of $y(x)$ by performing a random walk in which $+1$ corresponds to a pair of cousin primes and $-1$ to a pair of twin primes encountered by the walker. In particular, Wolf observed that $y(x)=0$ for $2823290$ primes up to $x=2^{43}\approx8.8\times10^{12}$ and noted that the clusters of $x$ where $y(x)=0$ formed a self-similar set.
	
	\subsection{Scope of this paper}
	
	The aim of this paper is threefold; first, to compute the Skewes number for the prime tuples listed below in Table \ref{table-tuples}, second, to provide a flexible algorithm for studying and gathering data on the asymptotic density of prime $k$-tuples and third, to provide computational data in support of the first Hardy-Littlewood conjecture. 
	
	In the remainder of this paper we shall denote by Li$_k(n)$ the integral $\displaystyle \int_{2}^{n} \frac{dt}{\log^{k}t}$ and by $\pi_P(n)$ the prime-counting function related to the prime tuple $P$. Table \ref{table-tuples} shows the admissible prime $k$-tuples we consider in this paper along with their conjectured asymptotic densities. Note that we include the twin primes within our study for the sake of completeness and in order to confirm Wolf's result that the corresponding Skewes number is indeed $1369391$.
	
	\begin{table}[h!] \caption{Prime $k$-tuples considered in this paper} \label{table-tuples} 
		\begin{tabular}{|c|l|l|}
			\hline
			Tuple & Definition & Conjectured asymptotic density \\ \hline
			$P_{2a}$ & $(p,p+2)$ & $\displaystyle 2 \prod_{p\geq3} \frac{p(p-2)}{(p-1)^2}$Li$_2(n)$ \\ \hline
			$P_{2b}$ & $(p,p+4)$ & $\displaystyle 2 \prod_{p\geq3} \frac{p(p-2)}{(p-1)^2}$Li$_2(n)$ \\ \hline
			$P_{3a}$ & $(p,p+2,p+6)$ & $\displaystyle \frac{9}{2} \prod_{p\geq5} \frac{p^2(p-3)}{(p-1)^3}$Li$_3(n)$ \\ \hline
			$P_{3b}$ & $(p,p+4,p+6)$ & $\displaystyle \frac{9}{2} \prod_{p\geq5} \frac{p^2(p-3)}{(p-1)^3}$Li$_3(n)$ \\ \hline
			$P_{4a}$ & $(p,p+2,p+6,p+8)$ & $\displaystyle \frac{27}{2} \prod_{p\geq5} \frac{p^3(p-4)}{(p-1)^4}$Li$_4(n)$\\ \hline
			$P_{4b}$ & $(p,p+4,p+6,p+10)$ & $\displaystyle 27 \prod_{p\geq5} \frac{p^3(p-4)}{(p-1)^4}$Li$_4(n)$ \\ \hline
			$P_{5a}$ & $(p,p+2,p+6,p+8,p+12)$ & $\displaystyle \frac{15^4}{2^{11}} \prod_{p\geq7} \frac{p^4(p-5)}{(p-1)^5}$Li$_5(n)$ \\ \hline
			$P_{5b}$ & $(p,p+4,p+6,p+10,p+12)$ & $\displaystyle \frac{15^4}{2^{11}} \prod_{p\geq7} \frac{p^4(p-5)}{(p-1)^5}$Li$_5(n)$ \\ \hline
			$P_{6}$  & $(p,p+4,p+6,p+10,p+12,p+16)$ & $\displaystyle \frac{15^5}{2^{13}} \prod_{p\geq7} \frac{p^5(p-6)}{(p-1)^6}$Li$_6(n)$ \\ \hline
		\end{tabular}
	\end{table}

	\section{Computational and algorithmic methods} \label{sec-methods}
	
	In this section we present the computational and algorithmic methods we used to obtain our results in Section \ref{sec-results}. These include finding the Skewes number for the prime $k$-tuples listed in Table \ref{table-tuples} as well as gathering other relevant numerical data in support of Conjecture \ref{conj-asym-dens}. 
	
	Note that the pseudo-code presented in this section, supplemented with code used for plotting purposes and gathering other secondary data was implemented in Wolfram Mathematica $11.1$ and executed on an Intel Core i7-$7800$X CPU @$3.50$ GHz and $16$ GB RAM to produce the results in Section \ref{sec-results}.
	
	The core of our algorithm relies on computing the logarithmic integral within successive intervals, the upper and lower bounds corresponding to consecutive primes of a given $k$-tuple. This method was also employed by Wolf \cite{Wolf11} to compute the Skewes number for twin primes. Unfortunately he did not give any code in support of his algorithm and did not consider any special cases that might arise such as insufficient computing memory or other extensions that allow recording additional data. Here we present two algorithms as pseudo-code. The first outlines our basic methodology while the second considers limitations in computing resources. Both of these algorithms are flexible, i.e., they can be easily modified to record intermediary data and fine-tune the parameters.
	
	We will now follow with a short presentation of the main algorithm. We begin by defining our inputs. These are:
	\begin{itemize}
		\item An array of primes in the prime $k$-tuple \textbf{P} in the range $[2,n]$,
		\item The corresponding prime tuple constant \textbf{C}, shown in Table \ref{table-tuples},
		\item The corresponding logarithmic integral function \textbf{LiP}$[]$.
	\end{itemize}
	The main program loop then computes the logarithmic integral within successive intervals and stores the cumulated value in a variable. This value is then multiplied by \textbf{C} and compared to the current value of $\pi_P$, which conveniently is the index of the current prime in \textbf{P}. The program continues until a reversal of the sign in the difference between $\pi_P$ and \textbf{C} times the cumulated logarithmic interval is reached. At this point, the Skewes number is returned. Algorithm \ref{algo-skewes} illustrates this method. Please note that the first prime is denoted by $p_1$ in the algorithm below.
	
	\begin{algorithm}[h!]
		\caption{Finding the Skewes number for a prime $k$-tuple}\label{algo-skewes}
		
		\begin{algorithmic}[0]
			
			\State \textbf{Input:} Prime $k$-tuple \textbf{P} within the range $[2,n]$, prime tuple constant \textbf{C}, \textbf{LiP}$\left[\ \right]$ function
			\State \textbf{Output:} Skewes number \textbf{S}
			
			\State $cumulLogIntegral \gets 0$
			\State $currentSign \gets $ Sign$[-1]$
			
			\For{\textbf{all} $p_i \in$ \textbf{P}}
			
			\If{$i=1$}
			\State $intLowerBound \gets 2$
			\Else
			\State $intLowerBound \gets p_{i-1}$
			\EndIf
			
			\State $intUpperBound \gets p_i$	
			\State $cumulLogIntegral \gets$ \textbf{LiP}[$intLowerBound$, $intUpperBound$]
			\State $difference \gets (i - \textbf{C} \times cumulLogIntegral)$
			
			\State $newSign \gets $ Sign$[difference]$
			
			\If{$currentSign \neq newSign$}
			\State \Return $p_i$
			\EndIf	
			
			\EndFor
			
		\end{algorithmic}
	\end{algorithm}
	
	A number of non-essential features were omitted from Algorithm \ref{algo-skewes}. For instance, it is desirable to record some of the intermediary values computed within the main loop such as the differences between $\pi_P(i)$ and its conjectured density. For instance, the last \textbf{if} clause can be enlarged with a block of code allowing to count the number of sign changes instead of returning the Skewes number. Our Mathematica implementation of the algorithm does just that, and some plots showing this data are presented in Section \ref{sec-results}.
	
	Obviously the biggest challenge in Algorithm \ref{algo-skewes} lies in providing the input to the algorithm, i.e., a list of primes within a given prime $k$-tuple up to a very high numerical limit. We note that the advantage of pre-computing such a list of primes instead of iterating through integers and testing for primality of other potential members of the tuple lies in the fact that highly efficient algorithms already exist within most mathematical software packages that are able to yield such an array in a computationally short amount of time. For instance, Wolfram Mathematica has the following one-line solution (taking the example of cousin primes):
	
	\begin{verbatim}
	
	cousinPrimes = Select[Prime[Range[PrimePi[upperLimit]]], PrimeQ[# + 4] &];
	
	\end{verbatim}
	with \texttt{upperLimit} replaced by any numerical upper bound. Of course, the average computer will quickly begin to struggle with the above code as larger values of \texttt{upperLimit} require generous amounts of memory to hold all the primes. For this reason, we present a modified version of Algorithm \ref{algo-skewes} such that, when a large enough upper bound is reached, we discard the contents of the current prime array and fill it with another ''chunk'' of equally large size; the current lower bound thus becomes the previous upper bound. We then use an offset in order to compensate for the prime counts in previous chunks. This revised version of the algorithm is presented in Algorithm \ref{algo-offset}.
	
	\begin{algorithm}[h!]
		\caption{Finding the Skewes number for a prime $k$-tuple, taking into account memory limitations}\label{algo-offset}
		
		\begin{algorithmic}[0]
			
			\State \textbf{Input:} \textbf{chunkSize}: the amount of primes within a chunk, prime tuple constant \textbf{C}, \textbf{LiP}$\left[\ \right]$ function
			\State \textbf{Output:} Skewes number \textbf{S}
			
			\State $offset \gets 0$
			
			\For{$chunk \gets 2$, \textbf{increment} \textbf{chunkSize},}
			
			\State \textbf{P} $\gets$ prime $k$-tuple  within the range $[chunk,chunk$+\textbf{chunkSize}$]$
			
			\State $cumulLogIntegral \gets 0$
			\State $currentSign \gets $ Sign$[-1]$
			
			\For{\textbf{all} $p_i \in$ \textbf{P}}
			
			\If{$i=1$}
			\State $intLowerBound \gets 2$
			\Else
			\State $intLowerBound \gets p_{i-1}$
			\EndIf
			
			\State $intUpperBound \gets p_i$	
			\State $cumulLogIntegral \gets$ \textbf{LiP}[$intLowerBound$, $intUpperBound$]	
			\State $difference \gets (i + offset - \textbf{C} \times cumulLogIntegral)$
			
			\State $newSign \gets $ Sign$[difference]$
			
			\If{$currentSign \neq newSign$}
			\State \Return $p_i$
			\EndIf	
			
			\EndFor
			
			$offset \gets offset+$ Length$[$\textbf{P}$]$
			\EndFor
			
		\end{algorithmic}
	\end{algorithm}
	
	We will now present the results obtained by implementing and running our algorithm in Wolfram Mathematica $11.1$.
	
	\section{Results} \label{sec-results}
	
	Using the methods described in Section \ref{sec-methods} we found Skewes numbers for $8$ new prime $k$-tuples. Table \ref{table-results-skewes} shows these numbers.
	
	\begin{table}[h!] \centering \caption{The Skewes numbers for the prime $k$-tuples considered in our study} \label{table-results-skewes}
		\begin{tabular}{|c|r|}
			\hline
			Prime $k$-tuple & Skewes number \\ \hline
			$P_{2a}$ & $1369391$             \\ \hline
			$P_{2b}$ & $5206837$              \\ \hline
			$P_{3a}$ & $87613571$              \\ \hline
			$P_{3b}$ & $337867$  \\ \hline
			$P_{4a}$ & $1172531$ \\ \hline
			$P_{4b}$ & $827929093$ \\ \hline
			$P_{5a}$ & $21432401$              \\ \hline
			$P_{5b}$ & $216646267$              \\ \hline
			$P_{6}$  & $251331775687$  \\ \hline
		\end{tabular}
	\end{table}
	
\begin{figure}[h!] 
		\caption{Value of $\displaystyle \pi_{P_{4a}}(n) - C_{P_{4a}} \int_2^n \frac{dt}{\log^{4}t}$ within the interval $[2,10^{8}]$}
		\label{P4-fig1}
		\centering
		\includegraphics[]{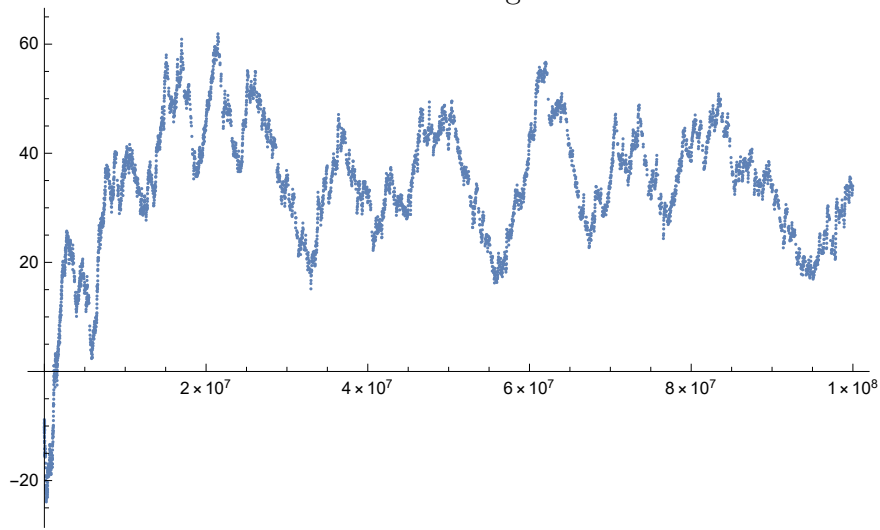}
	\end{figure}

\begin{figure}[h!]
		\caption{The crossover region in $P_{4a}$}
		\label{P4-fig2}
		\centering
		\includegraphics[]{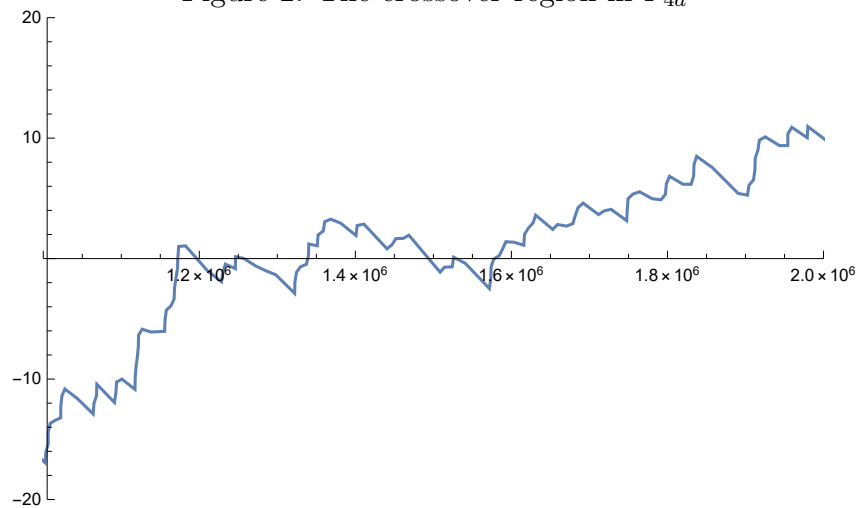}
	\end{figure}

	Various other results emerge from the above computations. First, we find that the inequality $\displaystyle \pi_P(n) < C_{a_1,a_2,...,a_k} \int_2^n \frac{dt}{\log^{k+1}t}$ only holds within a short interval for all of the prime tuples that have been considered. Taking for instance $P_{4a}$, it appears that the reversal of this inequality remains true after merely the $9^{\rm th}$ sign change of $\displaystyle \pi_{P_{4a}}(n) - C_{P_{4a}} \int_2^n \frac{dt}{\log^{4}t}$. Figure \ref{P4-fig1} shows a plot of this difference within the interval $[2,10^{8}]$ while Figure \ref{P4-fig2} zooms in on the crossover region. In another example, we observed a similar behaviour for $P_6$, which showed $15$ sign changes between its Skewes number ($251331775687$) and $26\times10^{10}$ (the Skewes number included).
	
	Our computations yielded other intriguing results as well. Indeed, it follows from the first Hardy-Littlewood conjecture that cousin primes ($P_{2b}=(p,p+4)$) have the same asymptotic density as twin primes, however the Skewes number for cousin primes is almost $4$ times greater than its twin prime counterpart. This is similar to the tuples $P_{5a}$ and $P_{5b}$ which, despite their equal conjectured density, have their Skewes number far apart; indeed, the one for $P_{5b}$ is almost $10$ times larger than that for $P_{5a}$. 
	
	Another consequence of the first Hardy-Littlewood conjecture is the thinning out of prime tuples within intervals of equal size. The results in our study do indeed provide evidence for such behaviour. Here we give the example of $P_6$ by comparing the value of $\pi_{P_{6}}(n)$ with the Hardy-Littlewood estimate (computed here to two significant figures after the decimal point) within intervals of size $10^{10}$. Table \ref{table-results-p6-intervals} shows our numerical results.
	
	\begin{table}[h!] \centering \label{table-results-p6-intervals} \caption{Values of $\pi_{P_{6}}(n)$ within intervals of size $10^{10}$ and the corresponding Hardy-Littlewood estimate (computed to two significant figures after the decimal point)}
		\begin{tabular}{|c|c|c|c|}
			\hline
			Interval ($[i,j]\times10^{10}$) & $\pi'_{P_{6}}(i,j)= \pi_{P_{6}}(j) - \pi_{P_{6}}(i)$ & Hardy-Littlewood & $\epsilon(i,j)/\pi'_{P_{6}}(i,j)$ \\ & & estimate $\epsilon(i,j)$ & \\ \hline
			$[0,1]$ & $1613$ &  $1664.55$   & $1.03195$\\ \hline
			$[1,2]$ & $1065$ & $1052.01$  & $0.98780$ \\ \hline
			$[2,3]$ & $897$  & $920.42$   & $1.02610$\\ \hline
			$[3,4]$ & $813$  & $845.55$   & $1.04003$\\ \hline
			$[4,5]$ & $796$  & $794.54$   & $0.99816$\\ \hline
			$[5,6]$ & $759$  & $756.47$   & $0.99666$\\ \hline
			$[6,7]$ & $674$  & $726.43$   & $1.07778$\\ \hline
			$[7,8]$ & $680$  & $701.82$   & $1.03208$\\ \hline
			$[8,9]$ & $680$  & $681.10$   & $1.00161$\\ \hline
			$[9,10]$ & $649$  & $663.29$  & $1.02201$\\ \hline
			$[10,11]$ & $638$  & $647.73$  & $1.01525$\\ \hline
			$[11,12]$ & $628$  & $633.96$  & $1.00949$\\ \hline
			$[12,13]$ & $610$  & $621.64$  & $1.01908$\\ \hline
			$[13,14]$ & $608$  & $610.52$  & $1.00414$\\ \hline
			$[14,15]$ & $605$  & $600.40$  & $0.99239$\\ \hline
			$[15,16]$ & $612$  & $591.13$  & $0.96589$\\ \hline
			$[16,17]$ & $613$  & $582.59$  & $0.95039$\\ \hline
			$[17,18]$ & $601$  & $574.69$  & $0.95622$\\ \hline
			$[18,19]$ & $620$  & $567.34$  & $0.91506$\\ \hline
		\end{tabular}
	\end{table}
	
	Finally, based on the results in our study, we propose the following conjecture.
	\begin{conjecture}
		All admissible prime $k$-tuples have a Skewes number.
	\end{conjecture}

	\section{Conclusion and further work}
	Alongside the discovery of $8$ new Skewes numbers for prime $k$-tuples, our results also give further evidence for the validity of the first Hardy-Littlewood conjecture. However, we find that tuples that are expected to have the same asymptotic density do not behave in the same manner, for instance their Skewes numbers occur considerable distances apart.
	
	Since our numerical data does not reach a sufficiently high upper bound, we are unable to investigate the validity of Wolf's conjecture within the context of other prime tuples. For instance, the cousin primes are expected to have the same density as the twin primes, thus according to Wolf's conjecture the amount of sign changes in the difference $\pi_{P_{2b}}(n) - C_{P_{2b}}$Li$_2(n)$ within the interval $n\in(1,T)$ should also be asymptotic to $\displaystyle \frac{\sqrt T}{\log T}$.
	
	We would also find interesting the extension of the search for Skewes numbers to further $k$-tuples such as $2$-tuples of the type $(p,p+2k)$ for positive integer $k>2$ (the case $k=3$ often being referred to as the ''sexy primes'' within more relaxed contexts), or longer ones such as
	\begin{itemize}
		\item $P_{7a} =$ $(p,p+2,p+6,p+8,p+12,p+18,p+20)$ and
		\item $P_{7b} =$ $(p,p+2,p+8,p+12,p+14,p+18,p+20)$. 
	\end{itemize}
	Both of these tuples are expected to have the same asymptotic density,
	$$
	\pi_{7a}(n),\pi_{7b}(n) \sim \displaystyle \frac{35^6}{3\times2^{22}} \prod_{p\geq11} \frac{p^6(p-7)}{(p-1)^7} \int_{2}^{n} \frac{dt}{\log^{7}t},
	$$
	but we did not find any primes that violate the Hardy-Littlewood inequality up to $1.2\times10^{11}$.
	
	\bibliographystyle{amsplain}

\begin{thebibliography}{9}
		\bibitem{Brent75}	
		R.~P.~Brent, \emph{Irregularities in the distribution of primes and twin primes}, Math. Comp. \textbf{29} (1975), 43--56.
		
		\bibitem{Wolf11}	
		M.~Wolf, \emph{The Skewes number for twin primes: counting sign changes of} $\pi_2(x)-C_2 \rm Li_2(x)$, Comput. Methods Sci. Technol. \textbf{17} (2011), 87--92.
		
		\bibitem{Wolf98}	
		M.~Wolf, \emph{Random walk on the prime numbers}, Physica A \textbf{250} (1998), 335--344.
		
		\bibitem{Nicely04}
		Th.~R.~Nicely, \emph{New evidence for the infinitude of some prime constellations}, 2004, \url{http://www.trnicely.net/ipc/ipc1d.pdf}
		
		\bibitem{HardyLittlewood22}
		G.~H.~Hardy, J.~E.~Littlewood, \emph{Some problems of 'Partitio Numerorum' III: On the expression of a number as a sum of primes}, Acta Math. \textbf{44}	(1922), 1--70.
		
		\bibitem{Littlewood14}
		J.~E.~Littlewood, \emph{Sur la distribution des nombres premiers}, C. R. Math. Acad. Sci. Paris \textbf{158} (1914), 1869--1872.
		
		\bibitem{Skewes33}
		S.~Skewes, \emph{On the difference} $\pi(x)-\rm {li} (x)$, J. London Math. Soc. \textbf{8} (1933), 277--283.
		
		\bibitem{Skewes55}
		S.~Skewes, \emph{On the difference} $\pi(x)-\rm {li} (x)$ (II), Proc. London Math. Soc. \textbf{5} (1955), 48--70.
		
		
	\end{thebibliography}

\end{document}